\documentclass[psamsfonts,reqno]{amsart}
\usepackage{amssymb,eucal,graphics,xy}
\usepackage{epstopdf}
\xyoption{all}

\theoremstyle{plain}

\newtheorem{lemma}{Lemma}[section]
\newtheorem{theorem}[lemma]{Theorem}

\newtheorem{corollary}[lemma]{Corollary}
\newtheorem{example}[lemma]{Example}

\theoremstyle{definition}

\theoremstyle{remark}

\newcommand{\qedsc}{{\qed}~{\rm Claim.}}

\numberwithin{equation}{section}

\newcommand{\set}[1]{\{#1\}}
\newcommand{\setm}[2]{\set{#1\mid#2}}

\newcommand{\famm}[2]{(#1\mid#2)}
\newcommand{\Famm}[2]{\left(#1\mid#2\right)}

\newcommand{\vin}{\mathbin{\vec{\in}}}
\newcommand{\Pow}{\mathfrak{P}}
\newcommand{\dnw}{\mathbin{\downarrow}}
\newcommand{\upw}{\mathbin{\uparrow}}

\newcommand{\cI}{\mathcal{I}}
\newcommand{\cF}{\mathcal{F}}
\newcommand{\into}{\hookrightarrow}
\newcommand{\onto}{\twoheadrightarrow}
\newcommand{\zero}[1]{{#1}^{\circ}}
\newcommand{\eps}{\varepsilon}
\newcommand{\go}{\omega}

\DeclareMathOperator{\card}{card}

\newcommand{\xp}{\mathbf{p}}
\newcommand{\xq}{\mathbf{q}}

\newcommand{\ol}[1]{\overline{{#1}}}

\newcommand{\res}{\mathbin{\restriction}}
\newcommand{\es}{\varnothing}
\newcommand{\op}{{\mathrm{op}}}

\DeclareMathOperator{\Eq}{Eq}
\DeclareMathOperator{\Id}{Id}
\DeclareMathOperator{\Fil}{Fil}
\newcommand{\Idb}{\operatorname{\ol{\Id}}}
\newcommand{\Filb}{\operatorname{\ol{\Fil}}}
\newcommand{\FL}{\mathrm{F}_{\mathbf{L}}}

\begin{document}

\title{Embedding coproducts of partition lattices}

\author[F.~Wehrung]{Friedrich Wehrung}
\address{LMNO, CNRS UMR 6139\\
D\'epartement de Math\'ematiques, BP 5186\\
Universit\'e de Caen, Campus 2\\
14032 Caen cedex\\
France}
\email{wehrung@math.unicaen.fr}
\urladdr{http://www.math.unicaen.fr/\~{}wehrung}


\dedicatory{Dedicated to B\'ela Cs\'ak\'any for his 75th birthday}

\date{\today}

\begin{abstract}
We prove that the lattice $\Eq\Omega$ of all equivalence relations on an infinite set~$\Omega$ contains, as a $0,1$-sublattice, the $0$-coproduct of two copies of itself, thus answering a question by G.\,M. Bergman. Hence, by using methods initiated by de Bruijn and further developed by Bergman, we obtain that $\Eq\Omega$ also contains, as a sublattice, the coproduct of $2^{\card\Omega}$ copies of itself.
\end{abstract}

\maketitle

\section{Introduction}\label{S:Intro}

Whitman's Theorem \cite{Whit} states that every lattice~$L$ can be embedded into the lattice $\Eq\Omega$ of all equivalence relations on some set~$\Omega$. The cardinality of~$\Omega$ may be taken equal to $\card L+\aleph_0$. There is not much room for improvement of the cardinality bound, as for example, $\Eq\Omega$ cannot be embedded into its dual lattice. (We believe the first printed occurrence of this result to be Proposition~6.2 in G.\,M. Bergman's recent preprint~\cite{Berg}, although it may have already been known for some time.) Hence the question of embeddability into~$\Eq\Omega$ of lattices of large cardinality (typically, $\card(\Eq\Omega)=2^{\card\Omega}$) is nontrivial.

In \cite{Berg}, Bergman also extends results of N.\,G. de Bruijn \cite{Brui1, Brui2} by proving various embedding results of large powers or copowers of structures such as symmetric groups, endomorphism rings, and monoids of self-maps of an infinite set~$\Omega$, into those same structures. The nature of the underlying general argument is categorical. The problem whether the lattice~$\Eq\Omega$ contains a coproduct (sometimes called ``free product'' by universal algebraists) of two, or more, copies of itself, was stated as an open question in a preprint version of that paper. In the present note, we solve this problem in the affirmative.

The idea of our proof is the following. The lattice~$\Eq\Omega$ of all equivalence relations on~$\Omega$ is naturally isomorphic to the ideal lattice~$\Id K$ of the lattice~$K$ of all \emph{finitely generated} equivalence relations, that is, those equivalence relations containing only finitely many non-diagonal pairs. Denote by~$K\amalg^0K$ the coproduct (amalgamation) of two copies of~$K$ above the common ideal~$0$. As $K\amalg^0K$ has the same cardinality as~$\Omega$, it follows from J\'onsson's proof of Whitman's Embedding Theorem that the lattice $\Id(K\amalg^0K)$ embeds into~$\Eq\Omega$. Finally, we prove that the ideal lattice functor preserves the coproduct~$\amalg^0$ and one-one-ness (Theorem~\ref{T:epsembedding}), in such a way that $(\Id K)\amalg^0(\Id K)$ embeds into~$\Id(K\amalg^0K)$. Then it is easy to extend this result to the usual coproduct $(\Id K)\amalg(\Id K)$.

We also present an example (Example~\ref{Ex:BamalgAC}) that shows that the result of Theorem~\ref{T:epsembedding} does not extend to amalgamation above a common (infinite) ideal. That is, for a common ideal~$A$ of lattices~$B$ and~$C$, the canonical homomorphism from $(\Id B)\amalg_{\Id A}(\Id C)$ to $\Id(B\amalg_AC)$ may not be one-to-one.

\section{Basic concepts}\label{S:Basic}

We refer to \cite{GLT2} for unexplained lattice-theoretical notions.
For any subsets~$Q$ and~$X$ in a poset (i.e., partially ordered set)~$P$, we put
 \[
 Q\dnw X=\setm{p\in Q}{(\exists x\in X)(p\leq x)}\quad\text{and}\quad
 Q\upw X=\setm{p\in Q}{(\exists x\in X)(p\geq x)}.
 \]
We also write $Q\dnw x$, resp. $Q\upw x$ in case $X=\set{x}$. A subset~$Q$ of~$P$ is a \emph{lower subset} of~$P$ if $Q=P\dnw Q$.

A map~$f\colon K\to L$ between lattices is \emph{meet-complete} if for each~$a\in K$ and each~$X\subseteq K$, $a=\bigwedge X$ in~$K$ implies that $f(a)=\bigwedge f[X]$ in~$L$. (Observe that we do not require either~$K$ or~$L$ to be a complete lattice.) When this is required only for nonempty~$X$, we say that~$f$ is \emph{nonempty-meet-complete}.
\emph{Join-completeness} and \emph{nonempty-join-completeness} of maps are defined dually. We say that~$f$ is \emph{complete} (resp., \emph{nonempty-complete}) if it is both meet-complete and join-complete (resp., both nonempty-meet-complete and nonempty-join-complete). We say that~$f$ is \emph{lower bounded} if $\setm{x\in K}{y\leq f(x)}$ is either empty or has a least element for each~$y\in L$. \emph{Upper bounded} homomorphisms are defined dually. Lower bounded homomorphisms are nonempty-meet-complete and upper bounded homomorphisms are nonempty-join-complete.

An \emph{ideal} of a lattice~$L$ is a nonempty lower subset of~$L$ closed under finite joins. We denote by~$\Id L$ the lattice of all ideals of~$L$. For a lattice homomorphism $f\colon K\to\nobreak L$, the map $\Id f\colon\Id K\to\Id L$ defined by
 \[
 (\Id f)(X)=L\dnw f[X]\,,\quad\text{for each }X\in\Id L\,,
 \]
is a nonempty-join-complete lattice embedding.
If $L$ is a \emph{$0$-lattice} (i.e., a lattice with least element), the canonical map $L\to\Id L$, $x\mapsto L\dnw x$ is a $0$-lattice embedding. The assignment that to every lattice associates its dual lattice~$L^\op$ (i.e., the lattice with the same underlying set as~$L$ but reverse ordering) is a category equivalence---and even a category isomorphism---from the category of all lattices to itself, that sends $0$-lattices to $1$-lattices. For every lattice~$L$, we denote by $\zero{L}$ the lattice obtained by adding a new zero element to~$L$.

A lattice $L$ is \emph{upper continuous} if for each $a\in L$ and each upward directed subset $\setm{x_i}{i\in I}$ of~$L$ admitting a join, the equality $a\wedge\bigvee_{i\in I}x_i=\bigvee_{i\in I}(a\wedge x_i)$ holds. We shall often use upper continuity in the following form: if~$I$ is an upward directed poset and both~$\famm{x_i}{i\in I}$ and $\famm{y_i}{i\in I}$ are isotone families with respective joins~$x$ and~$y$, then the family $\famm{x_i\wedge y_i}{i\in I}$ has join~$x\wedge y$.

Every algebraic lattice is upper continuous, so, for example, $\Id L\cup\set{\es}$ is upper continuous for any lattice~$L$; hence~$\Id L$ is also upper continuous. The lattice~$\Eq\Omega$ of all equivalence relations on a set~$\Omega$, partially ordered by inclusion, is an algebraic lattice, thus it is upper continuous. Other examples of upper continuous lattices that are not necessarily complete are given in~\cite{AGS}. For example, it follows from \cite[Corollary~2.2]{AGS} that \emph{every finitely presented lattice is upper continuous}.

We denote by $\Pow(\Omega)$ the powerset of a set~$\Omega$, and by~$\go$ the set of all natural numbers.

\section{The free lattice on a partial lattice}\label{S:FPLat}

We recall Dean's description of the free lattice on a partial lattice, see \cite{Dean} or \cite[Section~XI.9]{FJN}. A \emph{partial lattice} is a poset~$(P,\leq)$ endowed with partial functions~$\bigvee$ and~$\bigwedge$ from the nonempty finite subsets of~$P$ to~$P$ such that if $p=\bigvee X$ (resp., $p=\bigwedge X$), then~$p$ is the greatest lower bound (resp., least upper bound) of~$X$ in~$P$. An \emph{o-ideal} of~$P$ is a lower subset~$A$ of~$P$ such that $p=\bigvee X$ and $X\subseteq A$ implies that $p\in A$ for each $p\in P$ and each nonempty finite subset~$X$ of~$P$. The set~$\Idb P$ of all o-ideals of~$P$, partially ordered by inclusion, is an algebraic lattice. Observe that $\Idb P=(\Id P)\cup\set{\es}$ in case~$P$ is a lattice. O-filters are defined dually; again, the lattice~$\Filb P$ of all o-filters of~$P$, partially ordered by inclusion, is algebraic.
We denote by~$\cI(A)$ (resp., $\cF(A)$) the least o-ideal (resp., o-filter) of~$P$ containing a subset~$A$ of~$P$.

The \emph{free lattice}~$\FL(P)$ on~$P$ is generated, as a lattice, by an isomorphic copy of~$P$, that we shall identify with~$P$. (The subscript~$\mathbf{L}$ in~$\FL(P)$ stands for the variety of all lattices, as the ``free lattice on~$P$'' construction can be carried out in any variety of lattices.) For each~$x\in\FL(P)$, the following subsets of~$P$,
 \[
 \cI(x)=P\dnw x=\setm{p\in P}{p\leq x}\quad\text{and}
 \quad\cF(x)=P\upw x=\setm{p\in P}{x\leq p}
 \]
are, respectively, an o-ideal and an o-filter of~$P$, which can also be evaluated by the following rules:
 \begin{gather}
 \cI(x\vee y)=\cI(x)\vee\cI(y)\text{ in }\Idb P\,,\quad
 \cF(x\vee y)=\cF(x)\cap\cF(y);\label{Eq:cIcFxveey}\\
 \cI(x\wedge y)=\cI(x)\cap\cI(y),\quad\cF(x\wedge y)=
 \cF(x)\vee\cF(y)\text{ in }\Filb P\,,\label{Eq:cIcFxwedgey}
 \end{gather}
for all~$x,y\in\FL(P)$. The natural partial ordering on~$\FL(P)$ satisfies the following ``Whitman-type'' condition:
 \begin{multline}\label{Eq:Whit}
 x_0\wedge x_1\leq y_0\vee y_1\Longleftrightarrow\text{either }
 (\exists p\in P)(x_0\wedge x_1\leq p\leq y_0\vee y_1)\\
 \text{ or there is }i<2\text{ such that either }x_i\leq y_0\vee y_1\text{ or }x_0\wedge x_1\leq y_i\,,
 \end{multline}
which is also the basis of the inductive definition of that ordering.

\section{The $0$-coproduct of a family of lattices with zero}\label{S:0Coprod}

Our development of the~$0$-coproduct of a family of lattices with zero below bears some similarities with the development of coproducts (called there \emph{free products}) given in \cite[Chapter~VI]{GLT2}. Nevertheless, as we use the known results about the free lattice on a partial lattice (outlined in Section~\ref{S:FPLat}), our presentation becomes significantly shorter.

Let $\famm{L_i}{i\in I}$ be a family of lattices with zero. Modulo the harmless set-theoretical assumption that $L_i\cap L_j=\set{0}$ for all distinct indices~$i,j\in I$, the \emph{coproduct} (often called \emph{free product} by universal algebraists) of~$\famm{L_i}{i\in I}$ can be easily described as $\FL(P)$, where~$P$ is the partial lattice whose underlying set is the union~$\bigcup_{i\in I}L_i$, whose underlying partial ordering is the one generated by the partial orders on all the $L_i$s, and whose partial lattice structure consists of all existing joins and meets of nonempty finite subsets in each ``component''~$L_i$. We denote this lattice by $L=\coprod^0_{i\in I}L_i$, the superscript~$0$ meaning that the coproduct of the~$L_i$s is evaluated in the category of all $0$-lattices and $0$-preserving homomorphisms, which we shall often emphasize by saying ``$0$-coproduct'' instead of just coproduct. We shall also identify each~$L_i$ with its canonical copy in~$L$. Of course, the coproduct of any family of lattices $\famm{L_i}{i\in I}$ in the variety of all lattices is 
the sublattice of $\coprod^0_{i\in I}\zero{(L_i)}$ generated by the union of the images of the~$L_i$s.

Now we shall analyze further the structure of the $0$-coproduct~$L$, in a fashion similar to the development in \cite[Chapter~VI]{GLT2}. We add a new largest element, denoted by~$\infty$, to~$L$, and we set $\ol{L}_i=L_i\cup\set{\infty}$ for each~$i\in I$. The following lemma is an analogue, for $0$-coproducts instead of coproducts, of \cite[Theorem~VI.1.10]{GLT2}.

\begin{lemma}\label{L:xii}
For each~$x\in L$ and each~$i\in I$, there are a largest element of~$L_i$ below~$x$ and a least element of~$\ol{L}_i$ above~$x$ with respect to the ordering of~$L\cup\set{\infty}$. Furthermore, if we denote these elements by~$x_{(i)}$ and~$x^{(i)}$, respectively, then the following formulas hold:
 \begin{equation}\label{Eq:Defxii}
 \begin{aligned}
 &p_{(i)}=p^{(i)}=p,\text{ if }p\in L_i;\\
 &p_{(i)}=0\text{ and }p^{(i)}=\infty,\text{ if }p\in P\setminus L_i;\\
 &(x\vee y)_{(i)}=x_{(i)}\vee y_{(i)}
 \text{ and }(x\wedge y)_{(i)}=x_{(i)}\wedge y_{(i)};\\
 &(x\vee y)^{(i)}=x^{(i)}\vee y^{(i)};\\
 &(x\wedge y)^{(i)}=\begin{cases}
 0\,,&\text{if }x^{(j)}\wedge y^{(j)}=0\text{ for some }j\in I,\\
 x^{(i)}\wedge y^{(i)}\,,&\text{otherwise},
 \end{cases} 
 \end{aligned}
 \end{equation}
for each $x,y\in L$ and each $i\in I$.
\end{lemma}

\begin{proof}
For an element~$x$ of~$L$, abbreviate by ``$x_{(i)}$ exists'' (resp., ``$x^{(i)}$ exists'') the statement that $L_i\dnw x$ is a principal ideal in~$L_i$ (resp., $\ol{L}_i\upw x$ is a principal filter in~$\ol{L}_i$), and then denote by~$x_{(i)}$ (resp., $x^{(i)}$) the largest element of $L_i\dnw x$ (resp., the least element of~$\ol{L}_i\upw x$).
Denote by~$K$ the set of all $x\in L$ such that both $x_{(i)}$ and $x^{(i)}$ exist for each~$i\in I$. It is clear that~$K$ contains~$P$ and that both~$p_{(i)}$ and~$p^{(i)}$ are given by the first two formulas of~\eqref{Eq:Defxii}, for any $p\in P$. Furthermore, it follows immediately from the definition of~$K$ that
 \begin{align}
 \cI(z)&=\bigcup_{i\in I}(L_i\dnw z_{(i)}),\label{Eq:cI(x)fla}\\
 \cF(z)&=\bigcup_{i\in I}(L_i\upw z^{(i)}),\label{Eq:cF(x)fla}
 \end{align}
for each $z\in K$.
We shall establish that~$K$ is a sublattice of~$L$. So let $x,y\in K$, put $u=x\wedge y$ and $v=x\vee y$. It is straightforward that for each~$i\in I$, both~$u_{(i)}$ and $v^{(i)}$ exist, and
 \begin{equation}\label{Eq:u_iv^i}
 u_{(i)}=x_{(i)}\wedge y_{(i)},\quad v^{(i)}=x^{(i)}\vee y^{(i)}. 
 \end{equation}
Now we shall prove that $v_{(i)}$ exists and is equal to $x_{(i)}\vee y_{(i)}$.
By the induction hypothesis, \eqref{Eq:cI(x)fla} holds at both~$x$ and~$y$. So, as $\cI(v)=\cI(x)\vee\cI(y)$, in order to get the asserted existence and description of the elements~$v_{(i)}$, it suffices to prove that
 \begin{equation}\label{Eq:Lixiyivi}
 \bigcup_{i\in I}(L_i\dnw x_{(i)})\vee\bigcup_{i\in I}(L_i\dnw y_{(i)})=
 \bigcup_{i\in I}\bigl(L_i\dnw(x_{(i)}\vee y_{(i)})\bigr).
 \end{equation}
The containment from left to right is obvious, and each $x_{(i)}\vee y_{(i)}$ is contained in any o-ideal of~$P$ containing $\set{x_{(i)},y_{(i)}}$, so it suffices to prove that the right hand side of~\eqref{Eq:Lixiyivi} is an o-ideal of~$P$. As the join operation in~$P$ is internal to each~$L_i$, this set is closed under joins. As each~$L_i$ is a lower subset of~$P$, this set is also a lower subset of~$P$. This establishes the desired result for the~$v_{(i)}$s.

It remains to prove that $u^{(i)}$ exists and is equal to $z_i$, where $z_i=x^{(i)}\wedge y^{(i)}$ if $x^{(j)}\wedge y^{(j)}\neq0$ for all~$j$, and $z_i=0$ otherwise.
By the induction hypothesis, \eqref{Eq:cF(x)fla} holds at both~$x$ and~$y$. So, as $\cF(u)=\cF(x)\vee\cF(y)$, in order to get the asserted existence and description of the elements~$u^{(i)}$, it suffices to prove that
 \begin{equation}\label{Eq:Lixiyiui}
 \bigcup_{i\in I}(L_i\upw x^{(i)})\vee\bigcup_{i\in I}(L_i\upw y^{(i)})=
 \bigcup_{i\in I}(L_i\upw z_i).
 \end{equation}
The containment from left to right is obvious. If an o-filter~$U$ of~$P$ contains $\set{x^{(i)},y^{(i)}}$ for all~$i\in I$, then it also contains all elements $x^{(i)}\wedge y^{(i)}$; in particular, it is equal to~$P$ in case $x^{(i)}\wedge y^{(i)}=0$ for some~$i$. In any case, $z_i\in U$ for all~$i\in I$. So it suffices to prove that the right hand side of~\eqref{Eq:Lixiyiui} is an o-filter of~$P$. This is trivial in case~$z_i=0$ for some~$i$, so suppose that $z_i\neq0$ for all~$i$. As the meet operation in~$P$ is internal to each~$L_i$, the right hand side of~\eqref{Eq:Lixiyiui} is closed under meets. As each~$L_i\setminus\set{0}$ is an upper subset of~$P$, this set is also an upper subset of~$P$. This establishes the desired result for the~$u^{(i)}$s.
\end{proof}

\begin{lemma}\label{L:coprdKiLi}
Let~$K_i$ be a $0$-sublattice of a lattice~$L_i$, for each~$i\in I$. Then the canonical $0$-lattice homomorphism~$f\colon\coprod^0_{i\in I}K_i\to\coprod^0_{i\in I}L_i$ is an embedding.
\end{lemma}

\begin{proof}
By the amalgamation property for lattices \cite[Section~V.4]{GLT2}, the $i$-th coprojection from~$K_i$ to~$K$ is an embedding, for each~$i\in I$. Put $L'_i=K\amalg_{K_i}L_i$ for each~$i\in I$. Comparing the universal properties, it is immediate that the $0$-coproduct~$L$ of $\famm{L_i}{i\in I}$ is also the coproduct of $\famm{L'_i}{i\in I}$ over~$K$. Again by using the amalgamation property for lattices, all canonical maps from the~$L'_i$s to~$L$ are embeddings. So, in particular, the canonical map from their common sublattice~$K$ to~$L$ is an embedding.
\end{proof}

We shall call the adjoint maps $\alpha_i\colon x\mapsto x_{(i)}$ and $\beta_i\colon x\mapsto x^{(i)}$ the canonical lower, resp. upper adjoint of~$L$ onto~$L_i$, resp.~$\ol{L}_i$. Observe that these maps may not be defined in the case of amalgamation of two lattices over a common sublattice, as Example~\ref{Ex:BamalgAC} will show. (In that example, there is no largest element of~$B$ below $b_0\vee c_0$.)

The following result is an immediate consequence of well-known general properties of adjoint maps.

\begin{corollary}\label{C:Basicxii2}
The canonical embedding from $L_i$ into~$L$ is both lower bounded and upper bounded, for each~$i\in I$. In particular, it is a nonempty-complete lattice homomorphism. Furthermore, the lower adjoint $\alpha_i$ is meet-complete while the upper adjoint~$\beta_i$ is nonempty-join-complete.
\end{corollary}

In the following lemma, we shall represent the elements of $L=\coprod^0_{i\in I}L_i$ in the form~$\xp(\vec a)$, where~$\xp$ is a lattice term with variables from~$I\times\go$ and the ``vector'' $\vec a=\famm{a_{i,n}}{(i,n)\in I\times\go}$ is an element of the cartesian product~$\Pi=\prod_{(i,n)\in I\times\go}L_i$. Define a \emph{support} of~$\xp$ as a subset~$J$ of~$I$ such that~$\xp$ involves only variables from~$J\times\go$. Obviously, $\xp$ has a finite support. It is straightforward from~\eqref{Eq:Defxii} that $\xp(\vec a)_{(i)}=0$ and either~$\xp(\vec a)=0$ or $\xp(\vec a)^{(i)}=\infty$, for each~$i$ outside a support of~$\xp$.

\begin{lemma}\label{L:PartUC}
Let $\Lambda$ be an upward directed poset, let $\famm{{\vec a}^\lambda}{\lambda\in\Lambda}$ be an isotone family of elements of~$\Pi$ with supremum~$\vec a$ in~$\Pi$, and let~$\xp$ be a lattice term. If all the lattices~$L_i$ are upper continuous, then $\xp(\vec a)=\bigvee_{\lambda\in\Lambda}\xp({\vec a}^\lambda)$ in~$L$. 
\end{lemma}

Again, Example~\ref{Ex:BamalgAC} will show that Lemma~\ref{L:PartUC} fails to extend to the amalgam of two lattices over a common ideal.

\begin{proof}
As~$\xp(\vec a)$ is clearly an upper bound for all elements~$\xp({\vec a}^\lambda)$, it suffices to prove that for each lattice term~$\xq$ on~$I\times\go$ and each~$\vec b\in\Pi$ such that $\xp({\vec a}^\lambda)\leq\xq(\vec b)$ for all~$\lambda\in\Lambda$, the inequality~$\xp(\vec a)\leq\xq(\vec b)$ holds. We argue by induction on the sums of the lengths of~$\xp$ and~$\xq$. The case where~$\xp$ is a projection follows immediately from the second sentence of Corollary~\ref{C:Basicxii2}. The case where either~$\xp$ is a join or~$\xq$ is a meet is straightforward.

Now suppose that~$\xp=\xp_0\wedge\xp_1$ and~$\xq=\xq_0\vee\xq_1$. We shall make repeated uses of the following easily established principle, which uses only the assumption that~$\Lambda$ is upward directed:
 \begin{quote}\em
 For every positive integer~$n$ and every
 $X_0,\dots,X_{n-1}\subseteq\Lambda$, if
 $\bigcup_{i<n}X_i$ is cofinal in~$\Lambda$,
 then one of the~$X_i$s is cofinal in~$\Lambda$.
 \end{quote}
Now we use~\eqref{Eq:Whit}. If there exists a cofinal subset~$\Lambda'$ of~$\Lambda$ such that
 \[
 (\forall\lambda\in\Lambda')(\exists i<2)
 \bigr(\text{either }\xp_i({\vec a}^\lambda)\leq\xq(\vec b)\text{ or }
 \xp({\vec a}^\lambda)\leq\xq_i(\vec b)\bigr),
 \]
then there are~$i<2$ and a smaller cofinal subset~$\Lambda''$ of~$\Lambda'$ such that
 \begin{align*}
 \text{either}\quad&(\forall\lambda\in\Lambda'')\bigl(\xp_i({\vec a}^\lambda)\leq\xq(\vec b)\bigr)\\
 \text{or}\quad&(\forall\lambda\in\Lambda'')
 \bigl(\xp({\vec a}^\lambda)\leq\xq_i(\vec b)\bigr). 
 \end{align*}
In the first case, it follows from the induction hypothesis that~$\xp_i(\vec a)\leq\xq(\vec b)$. In the second case, it follows from the induction hypothesis that~$\xp(\vec a)\leq\xq_i(\vec b)$. In both cases, $\xp(\vec a)\leq\xq(\vec b)$.
It remains to consider the case where there exists a cofinal subset~$\Lambda'$ of~$\Lambda$ such that
 \[
 (\forall\lambda\in\Lambda')(\exists c_\lambda\in P)
 \bigl(\xp({\vec a}^\lambda)\leq c_\lambda\leq\xq(\vec b)\bigr).
 \]
It follows from the induction hypothesis that
 \begin{equation}\label{Eq:p01joinOK}
 \xp_\ell(\vec a)=\bigvee_{\lambda\in\Lambda'}\xp_\ell({\vec a}^\lambda)\,,
 \quad\text{for all }\ell<2\,.
 \end{equation}
Fix a common finite support~$J$ of~$\xp_0$, $\xp_1$, $\xq_0$, $\xq_1$.
Each~$c_\lambda$ belongs to~$L_i$, for some~$i$ in the given support~$J$. By using the finiteness of~$J$ and by extracting a further cofinal subset of~$\Lambda'$, we may assume that all those~$i$ are equal to the same index $j\in J$. Hence we have reduced the problem to the case where
 \begin{equation}\label{Eq:Interpolxpq}
 (\forall\lambda\in\Lambda')
 \bigl(\xp({\vec a}^\lambda)\leq c_\lambda\leq\xq(\vec b)\bigr)\,,
 \quad\text{where }c_\lambda=\xp({\vec a}^\lambda)^{(j)}\in L_j.
 \end{equation}
If $\xp_0(\vec a)^{(i)}\wedge\xp_1(\vec a)^{(i)}=0$ for some~$i\in I$, then~$\xp(\vec a)=0\leq\xq(\vec b)$ and we are done. Now suppose that $\xp_0(\vec a)^{(i)}\wedge\xp_1(\vec a)^{(i)}\neq0$ for all~$i\in I$. By using~\eqref{Eq:p01joinOK}, the finiteness of~$J$, and the upper continuity of~$L_i$, we obtain that there exists a cofinal subset~$\Lambda''$ of~$\Lambda'$ such that
 \[
 (\forall\lambda\in\Lambda'')(\forall i\in J)
 \bigl(\xp_0({\vec a}^\lambda)^{(i)}\wedge\xp_1({\vec a}^\lambda)^{(i)}
 \neq0\bigr).
 \]
In particular, both $\xp_0({\vec a}^\lambda)$ and $\xp_1({\vec a}^\lambda)$ are nonzero for each $\lambda\in\Lambda''$. As~$J$ is a common support of~$\xp_0$ and~$\xp_1$, the equality
$\xp_0({\vec a}^\lambda)^{(i)}\wedge\xp_1({\vec a}^\lambda)^{(i)}=\infty$ holds for all $\lambda\in\Lambda''$ and all $i\in I\setminus J$, hence
 \[
 (\forall\lambda\in\Lambda'')(\forall i\in I)
 \bigl(\xp_0({\vec a}^\lambda)^{(i)}\wedge\xp_1({\vec a}^\lambda)^{(i)}
 \neq0\bigr).
 \]
Thus it follows from~\eqref{Eq:Defxii} that $c_\lambda=\xp({\vec a}^\lambda)^{(j)}=\xp_0({\vec a}^\lambda)^{(j)}\wedge\xp_1({\vec a}^\lambda)^{(j)}$ for each $\lambda\in\Lambda''$. Hence, by the upper continuity of~$L_j$ (and thus of~$\ol{L}_j$), \eqref{Eq:p01joinOK}, and the previously observed fact that the upper adjoint~$\beta_j$ is nonempty-join-complete, $\setm{c_\lambda}{\lambda\in\Lambda''}$ has a join in~$L_j$, which is equal to~$\xp_0(\vec a)^{(j)}\wedge\xp_1(\vec a)^{(j)}=\xp(\vec a)^{(j)}$. Therefore, it follows from~\eqref{Eq:Interpolxpq} that $\xp(\vec a)\leq\xp(\vec a)^{(j)}\leq\xq(\vec b)$.
\end{proof}

\section{Ideal lattices and $0$-coproducts}\label{S:IdCoprod}

In this section we fix again a family~$\famm{L_i}{i\in I}$ of lattices with zero, pairwise intersecting in~$\set{0}$, and we form~$L=\coprod^0_{i\in I}L_i$. We denote by $\eps_i\colon\Id L_i\into\Id L$ the $0$-lattice homomorphism induced by the canonical embedding~$L_i\into L$, for each~$i\in I$. By the universal property of the coproduct, there exists a unique $0$-lattice homomorphism~$\eps\colon\coprod^0_{i\in I}\Id L_i\to\Id L$ such that $\eps_i=\eps\res_{\Id L_i}$ for each~$i\in I$. Observe that in case~$I$ is finite, the lattice~$\coprod^0_{i\in I}\Id L_i$ has $\bigvee_{i\in I}L_i$ as a largest element, and this element is sent by~$\eps$ to~$L$ (because every element of~$L$ lies below some join of elements of the~$L_i$s). Hence, \emph{if the index set~$I$ is finite, then the  map~$\eps$ preserves the unit as well}.

\begin{lemma}\label{L:Evalp(vecX)}
Let $\xp$ be a lattice term on~$I\times\go$ and let $\vec X=\famm{X_{i,n}}{(i,n)\in I\times\go}$ be an element of~$\prod_{(i,n)\in I\times\go}\Id L_i$. We put
$\vec\eps\vec X=\famm{\eps_i(X_{i,n})}{(i,n)\in I\times\go}\in(\Id L)^{I\times\go}$. Then the following equality holds.
 \[
 \xp(\vec\eps\vec X)=L\dnw\setm{\xp(\vec x)}{\vec x\vin\vec X}\,,
 \]
where ``\,$\vec x\vin\vec X$\,'' stands for $(\forall(i,n)\in I\times\go)(x_{i,n}\in X_{i,n})$.
\end{lemma}

\begin{proof}
We argue by induction on the length of the term~$\xp$. If~$\xp$ is a projection, then the result follows immediately from the definition of the maps~$\eps_i$. If~$\xp$ is either a join or a meet, then the result follows immediately from the expressions for the join and the meet in the ideal lattice of~$L$, in a fashion similar to the end of the proof of \cite[Lemma~I.4.8]{GLT2}.
\end{proof}

\begin{theorem}\label{T:epsembedding}
The canonical map~$\eps\colon\coprod^0_{i\in I}\Id L_i\to\Id\bigl(\coprod^0_{i\in I}L_i\bigr)$ is a $0$-lattice embedding.
\end{theorem}

\begin{proof}
We put again~$L=\coprod^0_{i\in I}L_i$.
Let~$\xp$, $\xq$ be lattice terms in~$I\times\go$ and let
$\vec X\in\prod_{(i,n)\in I\times\go}\Id L_i$ such that $\xp(\vec\eps\vec X)\leq\xq(\vec\eps\vec X)$ in~$\Id L$. We must prove that $\xp(\vec X)\leq\xq(\vec X)$ in
$\coprod^0_{i\in I}\Id L_i$. For each $\vec x\vin\vec X$, the inequalities $L\dnw\xp(\vec x)\leq\xp(\vec\eps\vec X)\leq\xq(\vec\eps\vec X)$ hold in~$\Id L$, thus, by Lemma~\ref{L:Evalp(vecX)}, there exists $\vec y\vin\vec X$ such that $L\dnw\xp(\vec x)\leq L\dnw\xq(\vec y)$ in~$\Id L$, that is, $\xp(\vec x)\leq\xq(\vec y)$ in~$L$. Therefore, by applying the canonical map from~$L=\coprod^0_{i\in I}L_i$ to~$\coprod^0_{i\in I}\Id L_i$ and putting $\vec L\dnw\vec x=\famm{L_i\dnw x_{i,n}}{(i,n)\in I\times\go}$, we obtain
 \begin{equation}\label{Eq:psmallqbig}
 \xp(\vec L\dnw\vec x)\leq\xq(\vec L\dnw\vec y)\leq\xq(\vec X)
 \quad\text{in }\coprod\nolimits^0_{i\in I}\Id L_i\,.
 \end{equation}
As $\vec X$ is equal to the directed join $\bigvee_{\vec x\vin\vec X}(\vec L\dnw\vec x)$ in~$\prod_{(i,n)\in I\times\go}\Id L_i$ and each $\Id L_i$ is upper continuous, it follows from Lemma~\ref{L:PartUC} that
 \[
 \xp(\vec X)=\bigvee\Famm{\xp(\vec L\dnw\vec x)}{\vec x\vin\vec X}\quad
 \text{in }\coprod\nolimits^0_{i\in I}\Id L_i\,.
 \]
Therefore, it follows from~\eqref{Eq:psmallqbig} that
 \begin{equation*}
 \xp(\vec X)\leq\xq(\vec X)\quad\text{in }\coprod\nolimits^0_{i\in I}\Id L_i\,.
 \tag*{\qed}
 \end{equation*}
 \renewcommand{\qed}{}
 \end{proof}
 
The following example shows that Theorem~\ref{T:epsembedding} does not extend to the amalgam $B\amalg_AC$ of two lattices~$B$ and~$C$ above a common ideal~$A$. The underlying idea can be traced back to Gr\"atzer and Schmidt in~\cite[Section~5]{GrSc95}.

\begin{example}\label{Ex:BamalgAC}
Lattices $B$ and $C$ with a common ideal~$A$ such that the canonical lattice homomorphism $f\colon(\Id B)\amalg_{\Id A}(\Id C)\to\Id\bigl(B\amalg_AC\bigr)$ is not one-to-one.
\end{example}

\begin{proof}
Denote by~$K$ the poset represented in Figure~\ref{Fig:LattK}. We claim that the subsets~$A$, $B$, and~$C$ of~$K$ defined by
 \begin{align*}
 A&=\setm{a_n}{n<\go}\cup\setm{p_n}{n<\go}\cup\setm{q_n}{n<\go}\,,\\
 B&=A\cup\setm{b_n}{n<\go}\,,\\
 C&=A\cup\setm{c_n}{n<\go}\,.
 \end{align*}
are as required. Observe that~$B$ and~$C$ are isomorphic lattices and that~$A$ is an ideal of both~$B$ and~$C$.

\begin{figure}[htb]
\includegraphics{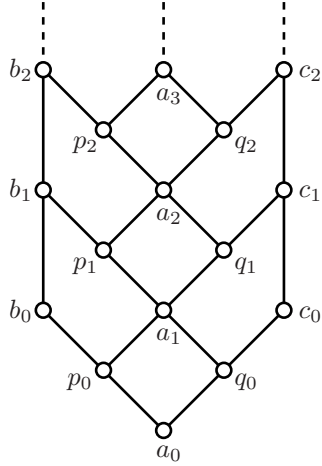}
\caption{The poset $K$.}\label{Fig:LattK}
\end{figure}

The map~$f$ is the unique lattice homomorphism that makes the diagram of Figure~\ref{Fig:CanMap} commute. Unlabeled arrows are the corresponding canonical maps.

\begin{figure}[htb]
 \[
 \xymatrix{
 & \Id B\ar[rd]\ar[rrrd] & & & \\
 \Id A\ar[ru]\ar[rd] & & (\Id B)\amalg_{\Id A}(\Id C) \ar[rr]|-{\textstyle f} & & \Id\bigl(B\amalg_AC\bigr)\\
 & \Id C\ar[ru]\ar[rrru] & & & \\
 }
 \]
\caption{The commutative diagram defining the homomorphism~$f$.}
\label{Fig:CanMap}
\end{figure}
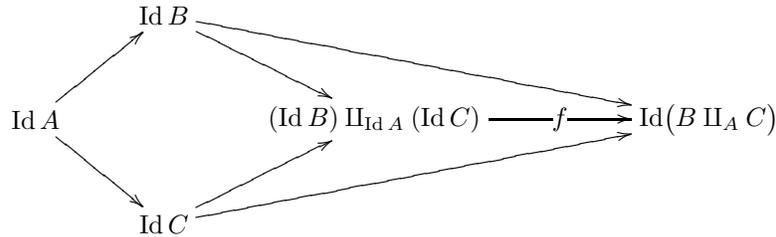

Put $D=B\amalg_AC$ and identify~$B$ and~$C$ with their images in~$D$. Further, we endow $\Id B\cup\Id C$ with its natural structure of partial lattice, that is, the ordering is the union of the orderings of~$\Id B$ and~$\Id C$ (remember that~$A=B\cap C$ is an ideal of both~$B$ and~$C$) and the joins and meets are those taking place in either~$\Id B$ or~$\Id C$. Observe that $\Id A=\Id B\cap\Id C$ and $(\Id B)\amalg_{\Id A}(\Id C)$ is the free lattice on the partial lattice~$(\Id B)\cup(\Id C)$. As the latter is identified with its canonical image in $(\Id B)\amalg_{\Id A}(\Id C)$, the elements~$A$, $B\dnw b_0$, and~$C\dnw c_0$ belong to $(\Id B)\amalg_{\Id A}(\Id C)$.

We prove by induction that $a_n\leq b_0\vee c_0$ in~$D$ for all~$n<\go$. This is trivial for $n=0$. Suppose that $a_n\leq b_0\vee c_0$. Then $a_n\vee b_0\leq b_0\vee c_0$, but~$B$ is a sublattice of~$D$ containing the subset~$\set{a_n,b_0}$ with join~$b_n$, thus $b_n\leq b_0\vee c_0$, and thus $p_n\leq b_0\vee c_0$. Similarly, $q_n\leq b_0\vee c_0$, but~$A$ is a sublattice of~$D$ containing the subset~$\set{p_n,q_n}$ with join~$a_{n+1}$, and thus $a_{n+1}\leq b_0\vee c_0$, which completes the induction step.

So we have established the inequality
 \begin{equation}\label{Eq:Aleqb+c}
 f(A)\leq f(B\dnw b_0)\vee f(C\dnw c_0)\quad\text{in }
 \Id\bigl(B\amalg_AC\bigr)=\Id D\,.
 \end{equation}
Now observe that $\setm{B\dnw x}{x\in B}\cup\setm{C\dnw y}{y\in C}$ is an o-ideal of the partial lattice $(\Id B)\cup(\Id C)$, containing~$\set{B\dnw b_0,C\dnw c_0}$ and to which~$A$ does not belong. Hence, $A\notin\cI(\set{B\dnw b_0,C\dnw c_0})$, which means that $A\nleq(B\dnw b_0)\vee(C\dnw c_0)$ in $(\Id B)\amalg_{\Id A}(\Id C)$. Therefore, by~\eqref{Eq:Aleqb+c}, $f$ is not an embedding.
\end{proof}

As observed before, this example shows that Lemma~\ref{L:PartUC} fails to extend to the amalgam of two lattices over a common ideal. Indeed, while $A=\bigvee_n(A\dnw a_n)$ in~$\Id B$, the same equality fails in $(\Id B)\amalg_{\Id A}(\Id C)$. The reason for this is that $A\dnw a_n\leq(B\dnw b_0)\vee(C\dnw c_0)$ for each~$n$, while $A\nleq(B\dnw b_0)\vee(C\dnw c_0)$.

\section{Embedding coproducts of infinite partition lattices}\label{S:EmbPart}

Whitman's Embedding Theorem states that every lattice embeds into~$\Eq\Omega$, for some set~$\Omega$. We shall use a proof of Whitman's Theorem due to B. J\'onsson \cite{Jons}, see also \cite[Section~IV.4]{GLT2}. The following result is proved there.

\begin{lemma}\label{L:Jonss}
For every lattice~$L$ with zero, there are an infinite set~$\Omega$ and a map $\delta\colon\Omega\times\Omega\to L$ satisfying the following properties:
\begin{enumerate}
\item $\delta(x,y)=0$ if{f} $x=y$, for all $x,y\in\Omega$.

\item $\delta(x,y)=\delta(y,x)$, for all $x,y\in\Omega$.

\item $\delta(x,z)\leq\delta(x,y)\vee\delta(y,z)$, for all $x,y,z\in L$.

\item For all $x,y\in\Omega$ and all $a,b\in L$ such that $\delta(x,y)\leq a\vee b$, there are $z_1,z_2,z_3\in\nobreak\Omega$ such that $\delta(x,z_1)=a$, $\delta(z_1,z_2)=b$, $\delta(z_2,z_3)=a$, and $\delta(z_3,y)=b$.
\end{enumerate}
\end{lemma}

Observe, in particular, that the map~$\delta$ is \emph{surjective}. Furthermore, a straightforward L\"owenheim-Skolem type argument (``keeping only the necessary elements in~$\Omega$'') shows that one may take $\card\Omega=\card L+\aleph_0$.

The following is the basis for J\'onsson's proof of Whitman's Embedding Theorem.

\begin{corollary}\label{C:Jonss}
For every lattice~$L$ with zero and every set~$\Omega$ such that $\card \Omega=\card L+\aleph_0$, there exists a complete lattice embedding from~$\Id L$ into~$\Eq\Omega$.
\end{corollary}

\begin{proof}
Any map~$\delta$ as in Lemma~\ref{L:Jonss} gives rise to a map $\varphi\colon\Id L\to\Eq\Omega$ defined by the rule
 \begin{equation}\label{Eq:Defvarphi}
 \varphi(A)=\setm{(x,y)\in\Omega\times\Omega}{\delta(x,y)\in A}\,,
 \quad\text{for each }A\in\Id L\,,
 \end{equation}
and conditions (1)--(4) above imply that~$\varphi$ is a complete lattice embedding.
\end{proof}

\begin{theorem}\label{T:Copr0Part}
Let~$\Omega$ be an infinite set. Then there exists a $0,1$-lattice embedding from $(\Eq\Omega)\amalg^0(\Eq\Omega)$ into~$\Eq\Omega$.
\end{theorem}

\begin{proof}
Denote by~$K$ the sublattice of~$\Eq\Omega$ consisting of all \emph{compact} equivalence relations of~$\Omega$. Thus the elements of~$K$ are exactly the equivalence relations containing only finitely many non-diagonal pairs. In particular, $\Eq\Omega$ is canonically isomorphic to~$\Id K$.

Now we apply Corollary~\ref{C:Jonss} to $L=K\amalg^0K$. As $\card L=\card\Omega$, we obtain a complete lattice embedding $\varphi\colon\Id L\into\Eq\Omega$. However, $\Id L=\Id(K\amalg^0K)$ contains, by Theorem~\ref{T:epsembedding} and the last sentence of the first paragraph of Section~\ref{S:IdCoprod}, a $0,1$-sublattice isomorphic to~$(\Id K)\amalg^0(\Id K)$, thus to~$(\Eq\Omega)\amalg^0(\Eq\Omega)$.
\end{proof}

For any nonempty set~$\Omega$, form $\ol{\Omega}=\Omega\cup\set{\infty}$ for an outside point~$\infty$. As there exists a retraction~$\rho\colon\ol{\Omega}\onto\Omega$ (pick $p\in\Omega$ and send~$\infty$ to~$p$), we can form a meet-complete, nonempty-join-complete lattice embedding $\eta\colon\Eq\Omega\into\Eq\ol{\Omega}$ by setting
 \[
 \eta(\theta)=\setm{(x,y)\in\ol{\Omega}\times\ol{\Omega}}
 {(\rho(x),\rho(y))\in\theta}\,,
 \quad\text{for each }\theta\in\Eq\Omega\,,
 \]
and $\eta$ sends the zero element of~$\Eq\Omega$ to a nonzero element of~$\Eq\ol{\Omega}$. Hence, in case~$\Omega$ is infinite, $\zero{(\Eq\Omega)}$ completely embeds into~$\Eq\Omega$. As $(\Eq\Omega)\amalg(\Eq\Omega)$ is the sublattice of $\zero{(\Eq\Omega)}\amalg^0\zero{(\Eq\Omega)}$ generated by the union of the images of~$\Eq\Omega$ under the two canonical coprojections, it follows from Theorem~\ref{T:Copr0Part} and Lemma~\ref{L:coprdKiLi} that $(\Eq\Omega)\amalg(\Eq\Omega)$ has a $1$-lattice embedding into~$\Eq\Omega$. If we denote by~$\theta$ the image of zero under this embedding, then $(\Eq\Omega)\amalg(\Eq\Omega)$ has a $0,1$-lattice embedding into~$\Eq(\Omega/{\theta})$, and thus, as $\card(\Omega/{\theta})\leq\card\Omega$, into~$\Eq\Omega$. Hence we obtain

\begin{theorem}\label{T:CoprPart}
Let~$\Omega$ be an infinite set. Then there exists a $0,1$-lattice embedding from $(\Eq\Omega)\amalg(\Eq\Omega)$ into~$\Eq\Omega$.
\end{theorem}

By applying the category equivalence~$L\mapsto L^\op$ to Theorems~\ref{T:Copr0Part} and~\ref{T:CoprPart} and denoting by~$\amalg^1$ the coproduct of $1$-lattices, we obtain the following result.

\begin{theorem}\label{T:CoprPartop}
Let~$\Omega$ be an infinite set. Then there are $0,1$-lattice embeddings from $(\Eq\Omega)^\op\amalg^1(\Eq\Omega)^\op$ into~$(\Eq\Omega)^\op$ and from $(\Eq\Omega)^\op\amalg(\Eq\Omega)^\op$ into~$(\Eq\Omega)^\op$.
\end{theorem}

By using the results of \cite{Berg}, we can now fit the copower of the optimal number of copies of~$L=\Eq\Omega$ into itself. The variety~$\mathbf{V}$ to which we apply those results is, of course, the variety of all lattices with zero.
The functor to be considered sends every set~$I$ to~$F(I)=\coprod^0_IL$, the $0$-coproduct of~$I$ copies of~$L$. If we denote by $e_i^I\colon L\into F(I)$ the $i$-th coprojection, then, for any map $f\colon I\to J$, $F(f)$ is the unique $0$-lattice homomorphism from~$F(I)$ to~$F(J)$ such that $F(f)\circ e_i^I=e_{f(i)}^J$ for all~$i\in I$. Observe that even in case both~$I$ and~$J$ are finite, $F(f)$ does not preserve the unit unless~$f$ is surjective. The condition labeled~(9) in \cite[Section~3]{Berg}, stating that every element of~$F(I)$ belongs to the range of~$F(a)$ for some $a\colon n\to I$, for some positive integer~$n$, is obviously satisfied. Hence, by \cite[Theorem~3.1]{Berg}, $F(\Pow(\Omega))$ has a $0$-lattice embedding into $F(\go)^\Omega$. Furthermore, it follows from \cite[Lemma~3.3]{Berg} that~$F(\go)$ has a $0$-lattice embedding into $\prod_{1\leq n<\go}F(n)$. By Lemma~\ref{L:coprdKiLi} and Theorem~\ref{T:Copr0Part}, each~$F(n)$ has a $0,1$-lattice embedding into~$L$. As, by the final paragraph of \cite[Section~2]{Berg}, $L^\Omega$ has a $0$-lattice embedding into~$L$, we obtain the following theorem.

\begin{theorem}\label{T:OptEmb}
Let $\Omega$ be an infinite set. Then the following statements hold:
\begin{enumerate}
\item $\coprod^0_{\Pow(\Omega)}\Eq\Omega$ has a $0$-lattice embedding into~$\Eq\Omega$.

\item $\coprod_{\Pow(\Omega)}\Eq\Omega$ has a lattice embedding into~$\Eq\Omega$.

\item $\coprod^1_{\Pow(\Omega)}(\Eq\Omega)^\op$ has a $1$-lattice embedding into~$(\Eq\Omega)^\op$.

\item $\coprod_{\Pow(\Omega)}(\Eq\Omega)^\op$ has a lattice embedding into~$(\Eq\Omega)^\op$.
\end{enumerate}
\end{theorem}

This raises the question whether $(\Eq\Omega)\amalg^1(\Eq\Omega)$ embeds into~$\Eq\Omega$, which the methods of the present paper do not seem to settle in any obvious way. More generally, we do not know whether, for a sublattice~$A$ of~$\Eq\Omega$, the amalgam $(\Eq\Omega)\amalg_A(\Eq\Omega)$ of two copies of~$\Eq\Omega$ over~$A$ embeds into~$\Eq\Omega$.

\section*{Acknowledgment}
I thank George Bergman for many comments and corrections about the successive versions of this note, which resulted in many improvements in both its form and substance.

\end{document}